\documentclass{amsart}
\usepackage[utf8]{inputenc}
\usepackage{amsmath,amsfonts,amsthm}
\usepackage[showonlyrefs]{mathtools}
\usepackage{tikz,tikzscale}
\usetikzlibrary{calc}
\usetikzlibrary{cd}
\usepackage{xcolor}
\usepackage{subcaption}
\usepackage[most]{tcolorbox}

\title[$C^1$-mapping]{A $C^1$-mapping based on finite elements on quadrilateral and hexahedral meshes}
\author{Daniel Arndt \and Guido Kanschat}

\let\phi\varphi
\newcommand{\rf}[1]{\widehat{#1}}
\def\mesh{\mathbb M}
\def\cell{K}
\def\refcell{\rf K}
\def\grad{\nabla}
\def\refgrad{\rf \nabla}
\def\ve{{\mathbf e}}
\def\vn{{\mathbf n}}
\def\vp{{\mathbf p}}
\def\vq{{\mathbf q}}
\def\vt{{\mathbf t}}
\def\vx{{\mathbf x}}
\def\vertex{{\mathbf v}}
\def\vertices{\mathcal V}
\def\node#1#2{\mathcal N_{#1}^{#2}}
\def\basis#1#2{\psi_{#1}^{#2}}
\def\refu{\rf u}
\def\refx{\rf x}
\def\refy{\rf y}
\def\refz{\rf z}
\def\vrefn{\rf{\mathbf n}}
\def\vreft{\rf{\mathbf t}}
\def\vrefx{\rf{\mathbf x}}

\theoremstyle{definition}
\newtheorem{algorithm}{Algorithm}

\begin{document}

\begin{abstract}
Finite elements of higher continuity, say conforming in $H^2$ instead of $H^1$, require a mapping from reference cells to mesh cells which is continuously differentiable across cell interfaces. In this article, we propose an algorithm to obtain such mappings given a topologically regular mesh in the standard format of vertex coordinates and a description of the boundary. A variant of the algorithm with orthogonal edges in each vertex is proposed. We introduce necessary modifications in the case of adaptive mesh refinement with nonconforming edges. Furthermore, we discuss efficient storage of the necessary data.
\end{abstract}

\maketitle

\section{Introduction}
The Bogner-Fox-Schmit element~\cite{BognerFoxSchmit1965} is usually referred to as an element on rectangles. While it is particularly simple there, already Petera and Pittman~\cite{PeteraPittman1994} presented an isoparametric version in 1994. Their mapping to the reference cell requires a global, smooth mapping from a domain with Cartesian cells to a deformed domain. Furthermore, they compute the derivative degrees of freedom by a global smoothing process. Here, we present a simple procedure generating such a mapping based on local information only involving cells adjacent to the same vertex. It extends to nonconforming meshes due to adaptive refinement.

It seems that smooth approximations and mappings have not attracted much attention in the finite element community until recently, most likely inspired by the advent of isogeometric methods. We refer in particular to Sangalli, Takacs, and their coworkers, discussing smooth piecewise spline approximation on irregular triangular and quadrilateral meshes in~\cite{CollinSangalliTakacs2016,KaplSangalliTakacs2017}. Since their analysis is based on Argyris-type degrees of freedom, they require at least polynomials of degree five. Here, we use degree three,  which becomes possible by a restriction to simpler topology.

An obvious purpose of a $C^1$-mapping is the $H^2$-conforming discretization with an isoparametric Bogner-Fox-Schmit element. This can be achieved in a fairly simple way by ensuring that the $2d$ edges leaving an interior vertex of a $d$-dimensional mesh are aligned pairwise.

When the mesh describes a smooth manifold in a higher dimensional space, the consistency of normal vector of adjacent cells is always an issue, which has originated fixes like averaging of these vectors. Such an averaging improves the consistency error, but does not produce a conforming method, while a $C^1$-mapping does.

More complicated is the case of the element introduced by Austin, Manteuffel, and McCormick~\cite{AustinManteuffelMcCormick04} on rectangles. We have introduced commuting interpolation operators for this element in~\cite{SharmaKanschat18}, which at the end of each edge involve the normal derivative of the normal vector component. On the other hand, due to the choice of the anisotropic polynomial space, the Hermitian/Lagrangian interpolation involved requires the tangential derivative along the orthogonal edge, which only works on rectangles, where the two edges are actually orthogonal. This property can be maintained if we impose additionally on the mapping that all edges leaving a vertex are either orthogonal or aligned. Therefore, below we describe two versions of the mapping, without and with orthogonalization.

The mapping without additional orthogonalization transfers easily to adaptive meshes with hanging nodes. There, position and derivative values are obtained by extension of the standard interpolation known from multilinear mapping. With additional orthogonalization, a $d$-cubic space for the mapping is not anymore sufficient and an enrichment is necessary.

The purpose of this article is twofold. On the one hand, we describe the algorithmic construction of a piecewise polynomial $C^1$-mapping and the conditions that have to be met. On the other hand, we discuss an efficient layout of data for such a mapping. In Section~\ref{sec:mapping}, we recall the notion of smooth mappings and then focus on the properties of polynomial mappings on a single cell. In Section~\ref{sec:global}, we discuss the construction of such mappings based on given vertex coordinates. We close with remarks on the implementation in Section~\ref{sec:implementation}.

\section{Mapped polynomials and degrees of freedom}
\label{sec:mapping}
We assume that the domain $\Omega \subset \mathbb R^d$ with $d=2,3$ is
subdivided into a mesh $\mesh$ of mesh cells $\cell$. While we will
discuss the exact shape of the cells below, we assume that the positions of the
vertices of the mesh and their connectivity are given. Furthermore, let the boundary
$\partial\Omega$ be smooth and its tangent plane be known in every point. Each mesh cell
$\cell$ will be represented by a smooth bijection
$\Phi_{\cell} : \refcell \to \cell$ from the reference cell
$\refcell = [-1,1]^d$ into the domain $\mathbb R^d$. The whole mesh is
then described by the function $\Phi$ on the Cartesian product of the
reference cell, such that for each mesh cell $\cell$ there is a factor
of this product with $\Phi_\cell(\refcell) = \cell$, thus
\begin{gather}
  \Phi\colon \prod_{K\in\mesh}\refcell \to \mathbb R^d,
\end{gather}
where the range of $\Phi$ approximates $\Omega$. Two cells $\cell_i$
and $\cell_j$ are called neighbors of each other, if the face
$F_{ij} = \cell_i \cap \cell_j$ has dimension $d-1$. We call the mesh
conforming, if each interior face $F_{ij}$ is a whole facet of each of
the adjacent cells. The discussion of nonconforming meshes is deferred
to Section~\ref{sec:hanging-nodes}. Let $P = [-2,2]^d$ and let $T$ map the
$2^d$-fold Cartesian product $\refcell\times\dots\times\refcell$ to
$P$ such that the reference cells $\refcell$ are mapped into one orthant of $P$ by translation by the vectors
$(\pm1,\dots,\pm1)^T$ with all possible combinations of signs and possibly rotation around their centers by multiples of $\pi/2$. A mesh
is called regular, if for every interior vertex $v$, there are $2^d$
mesh cells $K_{v_i}$ and a continuous mapping $\Phi_v$ such that
\begin{gather}
  \label{eq:regular-vertex}
  \begin{tikzcd}[ampersand replacement=\&]
    \refcell\times\dots\times\refcell
    \arrow[rr, "\Phi"]
    \arrow[rd, "T"]
    \&\& \bigcup K_{v_i}\\
    \& P
    \arrow[ru, "\Phi_{v}"']
  \end{tikzcd}.
\end{gather}
The set $P$ together with the mappings $\Phi_v$ forms an atlas for
$\Omega$. Therefore, we define that the mapping $\Phi$ is
$C^1$ if each of the mappings $\Phi_v$ is in $C^1(P)$.

\subsection{Mapped shape functions}
Shape functions on a single cell $\cell$ are defined by pull-back of
functions on $\refcell$, such that for $\vx=\Phi_\cell(\vrefx)$ there
holds
\begin{gather}
  \label{eq:5}
  \begin{split}
    u(\vx) &= \refu(\vrefx),\\
  \grad u(\vx)
  &= \refgrad\Phi_\cell(\vrefx)^{-T}
  \refgrad \refu(\vrefx)
  \end{split}
\end{gather}
where $(\refgrad\Phi_\cell)^{-T}$ is the transpose of the inverse of
$\refgrad\Phi_\cell$. For later use, we introduce the Jacobi
determinant $J_\cell = \operatorname{det} \refgrad\Phi_\cell$.
Furthermore, we will abbreviate $\Phi_{\cell}$ by $\Phi$ if no
confusion arises. We denote the components of $\Phi$ by $\phi_i$. Tangential and
normal vectors to the boundary of $\refcell$ are transformed to the
corresponding ones on $\cell$ by contravariant and covariant
transformation, respectively:
\begin{gather}
  \label{eq:6}
  \begin{split}
    \vt_K(\vx) &= \refgrad \Phi(\vrefx) \vreft_{\refcell} (\vrefx),\\
    \vn_K(\vx) &= \refgrad \Phi(\vrefx)^{-T} \vrefn_{\refcell} (\vrefx).
  \end{split}
\end{gather}


%

Continuity of function values can be achieved by standard techniques
of mapped finite elements. Here, we have to obtain continuity of first
derivatives by additional measures.  First investigating a vertex at
position $\vx$, we see that for every two cells $\cell_1$ and $\cell_2$
sharing a vertex, we obtain the additional condition
\begin{gather}
  \Bigl(\refgrad\Phi_{1}\bigl(\Phi_1^{-1}\vx\bigr)\Bigr)^{-T} \refgrad \refu_1(\vrefx)
  =
  \Bigl(\refgrad\Phi_{2}\bigl(\Phi_2^{-1}\vx\bigr)\Bigr)^{-T} \refgrad \refu_2(\vrefx).
\end{gather}
This condition can be used to establish conforming degrees of freedom on a mesh $\mesh_h$.
A simple choice is to require $\refgrad\Phi_{1}=\refgrad\Phi_{2}$ and thus decouple
$C^1$-continuity of the mapped element into $C^1$-continuity of the mapping and $C^1$-continuity of the
finite element with respect to the reference geometry.

For a face of codimension one shared by cell $\cell_1$ and $\cell_2$,
all tangential derivatives of $u_1$ and $u_2$ coincide just by
continuity of the function itself and the mapping. Thus, we have to ensure continuity of the normal
derivative along this face. By Taylor expansion, a necessary condition
is the continuity of $\partial_{\vt}\partial_{\vn} u(p)$ for any
tangential vector $\vt$ of the face. By continuity, this implies that $u$ must be $C^2$ at every vertex.

\subsection{Hermitian tensor product elements in arbitrary dimension}
For the standard $d$-linear mapping, the normal derivative of $\Phi$ at an interface is determined by the vertices not adjacent to this interface, namely vertices not shared by the two cells at the interface.
Therefore, continuity of this normal derivative can only be achieved using nonlocal information on a cell. Therefore, we use Hermitian interpolation and assign suitable values to the normal derivative in the vertices.
The choice of such values can be made differently, depending on the structure of the mesh.

The shape function space we are going to use with respect to the reference cell can be described as follows.
On the interval $[-1,1]$ with vertices $v_0 = -1$ and $v_1 = 1$, a cubic polynomial $p(\refx)$ is uniquely determined by the
values of the Hermite node functionals
\begin{gather}
\label{eq:1}
\begin{aligned}
  \node00(p) & = p(-1),& \node01(p) & = p'(-1),\\
  \node10(p) & = p(1), & \node11(p) & = p'(1),
\end{aligned}
\end{gather}
All degrees of freedom in one dimension are at the end of the interval
and the node functionals are dual to the well known basis functions
\begin{gather}
  \label{eq:2}
  \begin{aligned}
    4\basis00(\xi) & = \xi^3-3\xi+2,& 4\basis01(\xi) & = \xi^3-\xi^2-\xi+1,\\
    4\basis10(\xi) & = -\xi^3+3\xi+2, & 4\basis11(\xi) & = \xi^3+\xi^2-\xi+1.
  \end{aligned}
\end{gather}
Higher dimensional versions of this element are obtained by simple tensor products.
In this case, the degrees of freedom are all located in the vertices $\vertex_i \in \vertices$,
and $\vertices$ is the set of vertices of the reference cell $[-1,1]^d$.
In three dimensions for instance, we obtain in each vertex the node functionals
\begin{gather}
  \label{eq:3}
  \begin{aligned}
    \node i{000}(p) &= p(\vertex_i) \\
    \node i{100}(p) &= \partial_{\refx} p(\vertex_i) &
    \node i{010}(p) &= \partial_{\refy} p(\vertex_i) &
    \node i{001}(p) &= \partial_{\refz} p(\vertex_i) \\
    \node i{110}(p) &= \partial_{\refx}\partial_{\refy} p(\vertex_i) &
    \node i{101}(p) &= \partial_{\refx}\partial_{\refz} p(\vertex_i) &
    \node i{011}(p) &= \partial_{\refy}\partial_{\refz} p(\vertex_i) \\
    \node i{111}(p) &= \partial_{\refx}\partial_{\refy}\partial_{\refz} p(\vertex_i).
  \end{aligned}
\end{gather}

We refer to the number of ones in the upper index as the order of the
node functional that coincides with the order of the derivative
involved.  The associated basis functions for instance for
$v_0 = (-1,-1,-1)^T$ are
\begin{gather}
  \label{eq:tensor-basis}
  \basis0{\kappa\lambda\mu}(\refx,\refy,\refz)
  = \basis0\kappa(\refx)\basis0\lambda(\refy)\basis0\mu(\refz).
\end{gather}
This may suffice as an example.
The general case can easily be constructed using binary vertex indices.

In arbitrary dimension, we use a binary multi-index
$\kappa = (\kappa_1,\dots,\kappa_d)$ 
and obtain the representation
\begin{gather}
  \label{eq:4}
  p(\vrefx) = \sum_{\vertex_i\in \vertices} \sum_{\kappa\in \mathcal I}
  \alpha_i^\kappa \prod_j \basis{i}{\kappa_j}(\refx_j).
\end{gather}

\subsection{Basic properties of Hermitian-like mapping}

Consider mapping functions $\Phi_\cell: [-1,1]^d \to \cell$,
which are vector-valued Hermitian-like tensor product elements
themselves. That is, each component of the vector $\Phi_\cell(\vrefx)$
is of the form in equation~\eqref{eq:4}.

This choice has the immediate consequence, that the values of
$\Phi_\cell$ and its Jacobian on any face of the cell are determined
by degrees of freedom in vertices adjacent to this face only. Thus,
continuity of the Jacobian can be achieved by choosing them consistently on the cells sharing a vertex. Consider for
instance a horizontal edge of the reference square, bounded by
vertices say $\vertex_0 = (-1,-1)^T$ and $\vertex_1 = (1,-1)^T$. Then,
the position of a point $\vx = \Phi(\refx,-1)$ on the image of this
edge is a vector valued, cubic polynomial, determined by the four
degrees of freedom $\Phi(\pm1,-1)$ and
$\partial_{\refx}\Phi(\pm1,-1)$. The tangential derivative
$\partial_{\refx}\Phi(\refx,-1)$ is uniquely determined by the same
degrees of freedom. The normal derivative
$\partial_{\refy}\Phi(\refx,-1)$ is a cubic polynomial determined by
the degrees of freedom $\partial_{\refy}\Phi(\pm1,-1)$ and
$\partial_{\refx\refy}\Phi(\pm1,-1)$. Accordingly, $\nabla\Phi$ on
this edge has the following form
\begin{gather}
  \label{eq:8}
  \nabla\Phi(\xi,-1) =
  \bigl(
    \partial_{\refx} \vp(\xi), \vq(\xi)
  \bigr),
\end{gather}
with the vector valued polynomials
\begin{gather}
  \label{eq:9}
  \begin{split}
    \vp(\refx) &= \Phi(\vertex_0) \basis00(\refx)
    + \partial_{\refx}\Phi(\vertex_0) \basis01(\refx)
    + \Phi(\vertex_1) \basis10(\refx)
    + \partial_{\refx}\Phi(\vertex_1) \basis11(\refx),\\
    \vq(\refx) &=  \partial_{\refy} \Phi(\vertex_0) \basis00(\refx)
    + \partial_{\refx\refy} \Phi(\vertex_0) \basis01(\refx)
    + \partial_{\refy} \Phi(\vertex_1) \basis10(\refx)
    + \partial_{\refx\refy} \Phi(\vertex_1) \basis11(\refx).
  \end{split}
\end{gather}

It remains to choose the actual node values for this
transformation. Obviously, $\Phi(\vertex_i)$ must be the position of
the corresponding vertex of the actual cell. The degrees of freedom
involving derivatives on the other hand are not predetermined by a
standard mesh geometry, with exception of the domain boundary.  As a
consequence of~\eqref{eq:6}, first order degrees of freedom
$\partial_{j}\Phi(\vertex_i)$ determine up to a factor
the directions of the edges attached to this vertex. At this point, we
observe that the degrees of freedom of higher order and their basis
functions only distort the geometry interior to each mesh cell without
changing their boundaries, therefore, we already constrain mappings
$\Phi_\cell$ by
\begin{gather}
  \label{eq:10}
  \node i\alpha(\Phi_\cell) = 0 \quad \forall \left|\alpha\right| > 1.
\end{gather}
Thus, the polynomial $\vq$ in the gradient of $\Phi$ reduces to
\begin{gather}
  \label{eq:11}
  \vq(\refx) = \partial_{\refy} \Phi(\vertex_0) \basis00(\refx)
  + \partial_{\refy} \Phi(\vertex_1) \basis10(\refx).
\end{gather}
We also could have used the second derivatives to simplify $\vq$, for
instance, making it linear. This is impractical though, since these
degrees of freedom influence two edges at a time. Thus, the derivative
on one edge would depend on an opposing vertex.

\section{Globally differentiable mapping of mesh cells}
\label{sec:global}
The main challenge for the construction of $C^1$-conforming elements on arbitrary quadrilaterals and hexahedra consists in finding a $C^1$-conforming mapping $\Phi_\cell$ from the reference cell $\widehat\cell$ to the grid cell $\cell$, namely continuity of $\Phi$ and $\nabla\Phi$ across interfaces. In order to be able to use polynomial orders as low as $d$-cubic, we confine ourselves to the case of meshes with regular vertices in the sense of equation~\eqref{eq:regular-vertex} with possibly adaptive refinement. Irregular vertices require mapping with higher polynomial degrees, see for instance~\cite{KaplSangalliTakacs2017}.

We follow the isoparametric approach, choosing a $d$-cubic mapping defined by Hermitian-like interpolation. Clearly, the values of $\Phi$ in the vertices $\widehat{\vx}_i$ of $\widehat\cell$ must be the vertices $\vx_i$ of $T$.
Furthermore, continuity of the tangential derivatives of $\Phi$ along interfaces is implied by continuity of $\Phi$ itself. It remains to ensure continuity of the normal derivative.

\subsection{Shape regular, conforming meshes}
\begin{figure}[tp]
  \centering
  \begin{tikzpicture}
  	\coordinate (c) at (0,0);
    \coordinate (t) at (.3,1.5);
    \coordinate (b) at (-.8,-1.5);
    \coordinate (l) at (-2,.4);
    \coordinate (r) at (2,.7);
    \coordinate (d1) at ($($(l)!1cm!(r)$)-(l)$);
    \coordinate (d2) at ($($(b)!1cm!(t)$)-(b)$);
    \draw (b) -- (c) -- (t);
    \draw (l) -- (c) -- (r);
    \draw[fill=black] (l) circle (1mm);
    \draw[fill=black] (c) circle (1mm);
    \draw[fill=black] (r) circle (1mm);
    \draw[fill=black] (b) circle (1mm);
    \draw[fill=black] (t) circle (1mm);
    \draw[red,dotted] (l) -- (r);
    \draw[red,thick] ($(c)-(d1)$) -- ($(c)+(d1)$);
    \draw[red,dotted] (b) -- (t);
    \draw[red,thick] ($(c)-(d2)$) -- ($(c)+(d2)$);

 	\coordinate (c) at (5,0);
    \coordinate (t) at (5.3,1.5);
    \coordinate (b) at (4.2,-1.5);
    \coordinate (l) at (3,.4);
    \coordinate (r) at (7,.7);
    \coordinate (d3) at ($(d1)+(d2)$);
    \coordinate (d4) at ($(d1)-(d2)$);
    \coordinate (aux1) at ($(c)+(d3)$);
    \coordinate (aux2) at ($(c)+(d4)$);
    \draw (b) -- (c) -- (t);
    \draw (l) -- (c) -- (r);
    \draw[fill=black] (l) circle (1mm);
    \draw[fill=black] (c) circle (1mm);
    \draw[fill=black] (r) circle (1mm);
    \draw[fill=black] (b) circle (1mm);
    \draw[fill=black] (t) circle (1mm);
    \draw[red,dashed] ($(c)-(d1)$) -- ($(c)+(d1)$);
    \draw[red,dashed] ($(c)-(d2)$) -- ($(c)+(d2)$);
    \draw[red,dotted] (c) -- ($(c)!2cm!(aux1)$);
    \draw[red,dotted] (c) -- ($(c)!-2cm!(aux1)$);
    \draw[red,dotted] (c) -- ($(c)!2cm!(aux2)$);
    \draw[red,dotted] (c) -- ($(c)!-2cm!(aux2)$);
    \draw[red,thick] (c) -- ($(c)!1cm!45:(aux1)$);
    \draw[red,thick] (c) -- ($(c)!-1cm!45:(aux1)$);
    \draw[red,thick] (c) -- ($(c)!1cm!45:(aux2)$);
    \draw[red,thick] (c) -- ($(c)!-1cm!45:(aux2)$);
  \end{tikzpicture}
  \caption{Determining the direction of edges leaving a vertex. Simple averaging (left) and orthogonalized (right). Auxiliary lines dotted, final directions in red. Dashed lines on the right are final directions on the left.}
  \label{fig:averaging}
\end{figure}
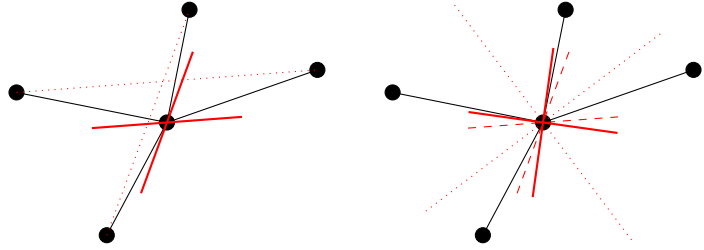

\begin{figure}[tp]
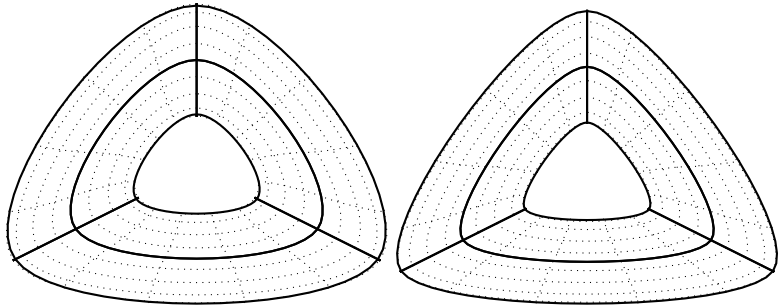

  \centering
  \includegraphics[width=.4\linewidth]{tikz/triangle_2.tikz}
  \includegraphics[width=.4\linewidth]{tikz/triangle_cubic.tikz}
  \caption{A triangular ring with bicubic mapping and different
    scaling of derivatives. Average distance from vertex to its
    neighbors (left) and distance of the two neighboring vertices
    (right).}
  \label{fig:triangle}
\end{figure}

Most easily, this can be achieved by choosing
the vector connecting the two adjacent vertices. In the situation of Figure~\ref{fig:averaging} on the left, the
direction vectors for the two edges leaving the center vertex are chosen parallel to the lines connecting the two vertices on these lines. This procedure is not local to each cell anymore, but requires averaging over the cells adjacent to a vertex.
The effect of such a mapping on a ring consisting of three trapezoidal cells can be seen in Figure~\ref{fig:triangle}. The shape regularity of the cell, indicated by distortions of the subdivision, is not much worse than for the corresponding bilinear mapping.

As can be seen in Figure~\ref{fig:averaging}, this method produces directions in a vertex which may not be orthogonal. While this is not a problem in general, it produces meshes not suitable for higher regularity divergence-conforming elements, see~\cite{AustinManteuffelMcCormick04,SharmaKanschat18}. Due to the anisotropic polynomial spaces used there, the normal derivative with respect to one edge in a corner point must be the tangent to the other edge. This requires that the edges meeting at one vertex are mutually orthogonal.

This can be achieved for instance by the Gram-Schmitt algorithm. Here, we propose the method illustrated in Figure~\ref{fig:averaging} on the right. First, we compute the lines dividing the angles between the two original directions in half. The directions chosen are obtained by rotating those by $\pi/2$. This way, both directions are close to the ones computed in the first step.

 After determining the direction, with or without orthogonalization, the norm of the derivative is chosen as the average length of the two edges it represents, as in Figure~\ref{fig:triangle} on the left. Even simpler, it can be chosen as half the distance of the two points used to compute the direction, as on the right. As the figure shows, simplicity comes at the cost of more strongly differing gradients.

\subsection{Hanging nodes}
\label{sec:hanging-nodes}
When dealing with locally refined meshes with nonconforming edges/faces, we follow the accepted approach that the degrees of freedom on the refined side must be constrained such that the mappings on both side coincide. This principle now applies to positions as well as derivatives.
\begin{figure}
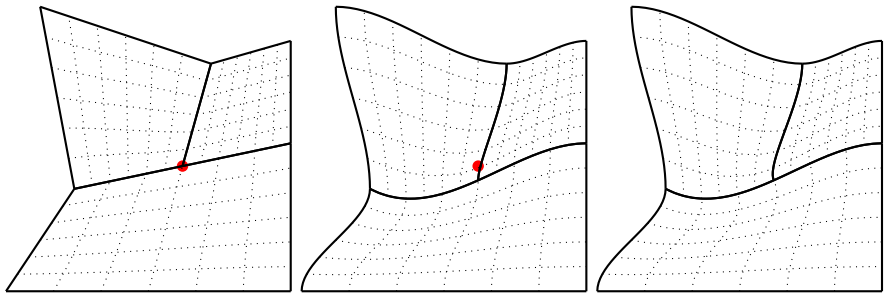

  \centering
  \includegraphics[width=.3\linewidth]{tikz/adaptive_linear.tikz}
  \includegraphics[width=.3\linewidth]{tikz/adaptive_cubic.tikz}
  \includegraphics[width=.3\linewidth]{tikz/adaptive_ortho.tikz}
  \caption{A one-irregular edge with bilinear $C^0$-mapping (left) and bicubic $C^1$-mapping (center) and augmented mapping with orthogonality at center point (right). The red dot indicates the location of the hanging vertex with linear interpolation on the edge.}
  \label{fig:one-irregular}
\end{figure}
Note, that the vertex in the center of the patch in Figure~\ref{fig:one-irregular} is artificial, hence the name \emph{hanging vertex}. Thus, its position can be chosen to match the requirements of the mapping.
While located at the average of its two neighbors' positions for standard mappings, its location is now determined by the cubic polynomial mapping on the coarse cell, see Figure~\ref{fig:one-irregular} left and center.
This approach has the advantage of being dimension independent, by taking the values of the $d-1$-cubic polynomials on the face of the coarse cell at the bisecting points. It can be easily extended to more levels of irregularity or to nonuniform subdivision.

The degrees of freedom of the mapping on one of the small cell in the hanging vertex are completely determined by the two vertices left and right of it. This includes not only the position and first derivatives, but also the mixed second derivative, which might be nonzero here.

We point out here, that the parameter $\xi$ on the large cell traverses the edge only once, while it does so twice on the two small cells. Therefore, the derivative degrees of freedom on the refined cells must be divided by two to maintain consistency with the coarse cell. Since adaptive algorithms typically generate a mesh hierarchy by refinement, this can be obtained by taking refinement level information into account. For details, refer to the data structures discussed in Section~\ref{sec:implementation}.

In order to be $C^1$-conforming, the direction of the edge going upward from the hanging vertex in Figure~\ref{fig:one-irregular} is determined by the normal derivative of the mapping on the coarse cell below, which in turn is uniquely determined by its values at the vertices. In particular, this edge may or may not be perpendicular to the other one. If we aim at vertices with orthogonal edges, we cannot avoid adding a shape function on the coarse side. Since this implies a 4th order polynomial along the edge, the same augmentation is required on the refined cells. Such an augmentation can be obtained by a hierarchical approach: first, the bicubic mapping on the large cell is computed according to the previous section. In a second step, degrees of freedom for additional shape functions are determined such that the vectors in $\refgrad \Phi$ are orthogonal in the centers of the edges (and the faces in 3D).

In two dimensions, we augment the mapping by one vector valued function for each edge, which allows adjusting the normal derivative. First, add to the one-dimensional basis in equation~\eqref{eq:2} the 4th order polynomial
\begin{gather}
	\psi_E(\xi) = 64 \prod_{i,j\in\{0,1\}} \psi_i^j(\xi),
\end{gather}
which has vanishing degrees of freedom $\node ij$ for $i,j\in\{0,1\}$ and the value 1 at $\xi=0$. Then, for each edge, add the mapping basis function, which is the product of $\psi_E$ in tangential direction with the basis function for the derivative degree of freedom in orthogonal direction. For example, we add the basis function
\begin{gather}
	\psi_{E,1}(\refx,\refy) = \psi_0^1(\refx) \psi_E(\refy),
\end{gather}
on the left edge $E_1 = \{-1\} \times [-1,1]$ which satisfies
\begin{align*}
 \psi_{E,1}(-1,\refy) &= 0, & \partial_{\refx} \psi_{E,1}(-1,\pm 1) &= 0, \\
 \partial_{\refy} \psi_{E,1}(-1,\refy) &= 0, & \partial_{\refx} \psi_{E,1}(-1,0) &= 1.
\end{align*}
For the derivatives at $x=-1, y=0$, this choice yields
\begin{align*}
 \partial_{\refx} \Phi &= \boldsymbol{\alpha}_0^{10}\partial_{\refx}\psi_0^{10}+\boldsymbol{\alpha}_0^{11}\partial_{\refx}\psi_0^{11}+\boldsymbol{\alpha}_2^{10}\partial_{\refx}\psi_2^{10}+\boldsymbol{\alpha}_2^{11}\partial_{\refx}\psi_2^{11}+\boldsymbol{\alpha}_{E,1}=:\boldsymbol{\alpha}_x +\boldsymbol{\alpha}_{E,1}\\
 \partial_{\refy} \Phi &= \boldsymbol{\alpha}_0^{01}\partial_{\refy}\psi_0^{01}+\boldsymbol{\alpha}_0^{11}\partial_{\refy}\psi_0^{11}+\boldsymbol{\alpha}_2^{01}\partial_{\refy}\psi_2^{01}+\boldsymbol{\alpha}_2^{11}\partial_{\refy}\psi_2^{11}=:\boldsymbol{\alpha}_y.
\end{align*}

A straightforward approach is to choose $\boldsymbol{\alpha}_{E,1}$ such that $\partial_{\refx} \Phi$ is the projection of $\boldsymbol{\alpha}_x$ onto a vector orthogonal to $\boldsymbol{\alpha}_y$.
Using the Gram-Schmidt procedure, this implies
\begin{gather}
\label{eq:Gram-Schmidt-1}
\partial_{\refx} \Phi = \boldsymbol{\alpha}_x - \frac{\boldsymbol{\alpha}_x\cdot\boldsymbol{\alpha}_y}{\|\boldsymbol{\alpha}_y\|^2}\boldsymbol{\alpha}_y
\quad\Leftrightarrow\quad
\boldsymbol{\alpha}_{E,1}   = - \frac{\boldsymbol{\alpha}_x\cdot\boldsymbol{\alpha}_y}{\|\boldsymbol{\alpha}_y\|^2}\boldsymbol{\alpha}_y.
\end{gather}

The construction extends by tensor product structure to three dimensional cells, where we have to consider face and edge augmentations of the form
\begin{gather}
	\psi_{E/F,\cdot}(\refx,\refy,\refz) = \psi_0^1(\refx) \psi_E(\refy) \psi_E(\refz),
\end{gather}
one on each face, two on each edge. The augmentation on the edges is chosen such that for each face the tangential vector not parallel to the edge is orthogonal to the edge. Then, these vectors form a two-dimensinoal system which can be orthogonalized according to Figure~\ref{fig:averaging}, obtaining an orthogonal system of three vectors, which form the directions of the edges of the refined cells.

After this ---note that we have modified the tangentials of the
faces--- the coefficient for the face augmentation
follows~\eqref{eq:Gram-Schmidt-1}, but orthogonalizing with respect to
both tangent vectors.

\subsection{Boundary layers}
\label{sec:anisotropic}
When discretizing problems exhibiting boundary layers, typically an
adjusted mesh with anisotropic mesh cells is used.
\begin{figure}[tp]
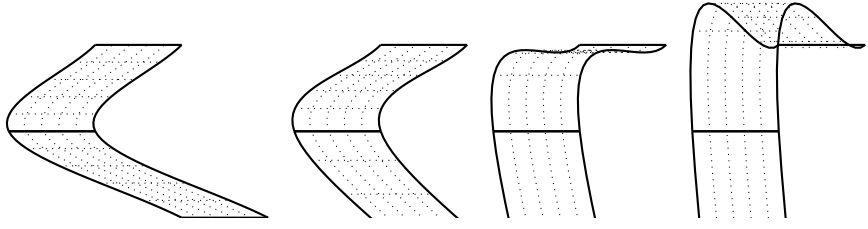

  \centering
  \includegraphics[height=.14\textheight]{tikz/angle1}
  \includegraphics[height=.14\textheight]{tikz/anglefail2}
  \includegraphics[height=.14\textheight]{tikz/anglefail8}
  \includegraphics[height=.14\textheight]{tikz/anglefail16}
  \caption[Example for failure with anisotropic cells]{An increasingly anisotropic wedge shaped mesh. The algorithm of the previous section is applied to the two center vertices. The bottom cells have heights 1, 2, 8, and 16, respectively (from left to right). Horizontal positions are constant.}
  \label{fig:anisotropic-1}
\end{figure}
The procedure described in the previous section averages the gradients of all
cells adjacent to a vertex. Thus, it can be expected that the approach results in self-intersecting cells if the mesh is not locally quasi-uniform, for instance, if a cell has a neighbor with high aspect ratio. This effect can indeed be observed in
Figure~\ref{fig:anisotropic-1}, where the meshes have been elongated to show the behavior of the mapping on the small cell more clearly. The original mesh on the left consists
of two cells of equal height, forming a wedge to highlight the effect
discussed. In the subsequent meshes, the lower cell is stretched to
towards the bottom more and more. While stretching by a factor 2
hardly affects the quality of the mapping in the smaller cell, the
mapping gets close to singular at a factor 8, where the
self-intersection could still be avoided by changing the derivatives
at the top edge. At a stretching of 16, the mapping is clearly
self-intersecting.

An obvious fix would be a scaling of the derivative of the mapping to
accommodate for the small cell. But, such a change would scale the
vertical derivative in the large cell, thus introducing a variation of
the Jacobi determinant $J$ depending on the small cell, which would in
turn yield bad constants in local interpolation estimates.

The solution is the choice of a mapping, which is not continuously
differentiable, but still yields such finite element functions by
transformation. Indeed, we do a similar rescaling of derivatives as in
the case of hanging vertices, but this time for the normal derivative
instead of the tangential one.

\begin{figure}[tp]
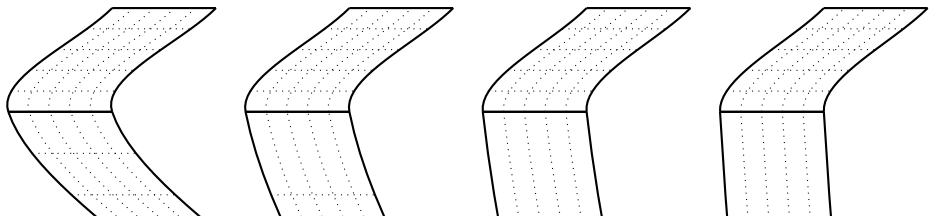

  \centering
  \includegraphics[width=.24\textwidth]{tikz/angle2}
  \includegraphics[width=.24\textwidth]{tikz/angle4}
  \includegraphics[width=.24\textwidth]{tikz/angle8}
  \includegraphics[width=.24\textwidth]{tikz/angle16}
  \caption[Example for anisotropic cells]{An increasingly anisotropic wedge shaped mesh with corrected mapping. The bottom cells have heights 2, 4, 8, and 16, respectively (from left to right). Horizontal positions are constant}
  \label{fig:anisotropic-2}
\end{figure}

This is achieved by defining on each of the two cells
$\cell_1$ and $\cell_2$ sharing a face $F$ a ``length orthogonal to
$F$'' called $h_{i;F}$. Several definitions of such a quantity can be
devised. Here, we suggest the following: for each vertex $\vertex$ of
$F$, compute the hyperplane obtained from $\vertex$ and all other
vertices of $F$ sharing an edge with $\vertex$. Then, compute the
distance of the (unique) vertex of $\cell_i$ sharing an edge with
$\vertex$ and not in $F$ to this hyperplane. Finally, take the average
over all vertices $\vertex$ of $F$. Figure~\ref{fig:anisotropic-2} shows cells with the same parameters as in Figure~\ref{fig:anisotropic-1}, but with rescaled derivatives on the small cell.
Clearly, the quality of the mapping does not deteriorate anymore when
the large cell is stretched.

The implementation of such a mapping can be achieved by an additional vector-valued scaling factor
stored for each cell and multiplied with the degree of freedom for the derivative.
However, note that the scaling factor in one direction must be shared by all cells joined by an edge in that direction.
Therefore, this scheme is suitable for boundary layers, where we have this situation, but not for anisotropies that may change side along a line.

\section{Implementation}
\label{sec:implementation}
The $C^1$-conforming meshes described above are implemented in the
open-source finite element library \texttt{deal.II}
\cite{dealII90,BangerthHartmannKanschat2007}. In particular,
the mapping is implemented in terms of the \texttt{MappingFEField}
class that allows to describe the geometry by a finite element vector.

The mesh data structure in $d$ dimensions requires $d+1$ vector-valued
degrees of freedom in each non-hanging vertex, one of order zero for
the position of the vertex and $d$ of order one for the
derivatives. While Hermite interpolation uses even higher order degrees of freedom, these are set to zero for mapping purposes.
 Using the local numbering of degrees of freedom described
in (\ref{eq:3}), these are assigned by the following algorithm. Note
that the algorithm uses additional intermediate data in each vertex,
which can be deleted after completion in order to free space for
numerical computation.
\begin{algorithm}
  Let a mesh geometry be given by vertex coordinates as well as the
  tangent planes in boundary vertices.
\begin{enumerate}
\item Assign the coordinates of each vertex to the zero order degree
  of freedom.
\item For each edge between a vertex and a neighboring vertex, record
  the connecting vector $\ve$ and a scaling factor $h_E$ taking into account
  anisotropic boundary cells or locally refined cells. Ignore hanging
  vertices in this process.

  For instance, this can be achieved cell-wise and the first vertex
  $\boldsymbol{v}_0$ on a cell in two dimensions gets the direction vectors
\begin{align*}
\ve_0^0 &=\boldsymbol{v}_1-\boldsymbol{v}_0\\
\ve_0^1 &=\boldsymbol{v}_2-\boldsymbol{v}_0
\end{align*}
assuming a lexicographic ordering of the vertices of the cell.
Note that this choice yields the standard $d$-linear $C^0$-mapping with discontinuous derivatives.
\item For each vertex and each pair of edges $e_1$ and $e_2$ reaching
  the vertex from opposite sides, compute the corresponding first
  order degree of freedom as the weighted mean
  \begin{gather*}
    \boldsymbol{\alpha}_v^{*}
    = \frac{\frac1{h_1}\ve_1 + \frac1{h_2}\ve_2}{\frac1{h_1}+\frac1{h_2}}.
  \end{gather*}
  For a shape regular mesh without hanging vertices, this is simply
  the average.

  If the vertex is on the boundary, project the obtained vector to the
  tangent plane.
\item Optionally: orthogonalize the derivative degrees of freedom in
  each vertex. If the mesh has hanging vertices, compute the
  additional degrees of freedom according to equation~\eqref{eq:Gram-Schmidt-1}.
\item Compute the positions and derivatives for hanging vertices
  according to Section~\ref{sec:hanging-nodes}.
 \end{enumerate}
\end{algorithm}

\begin{figure}
\begin{subfigure}{.49\linewidth}
\includegraphics[width=\linewidth]{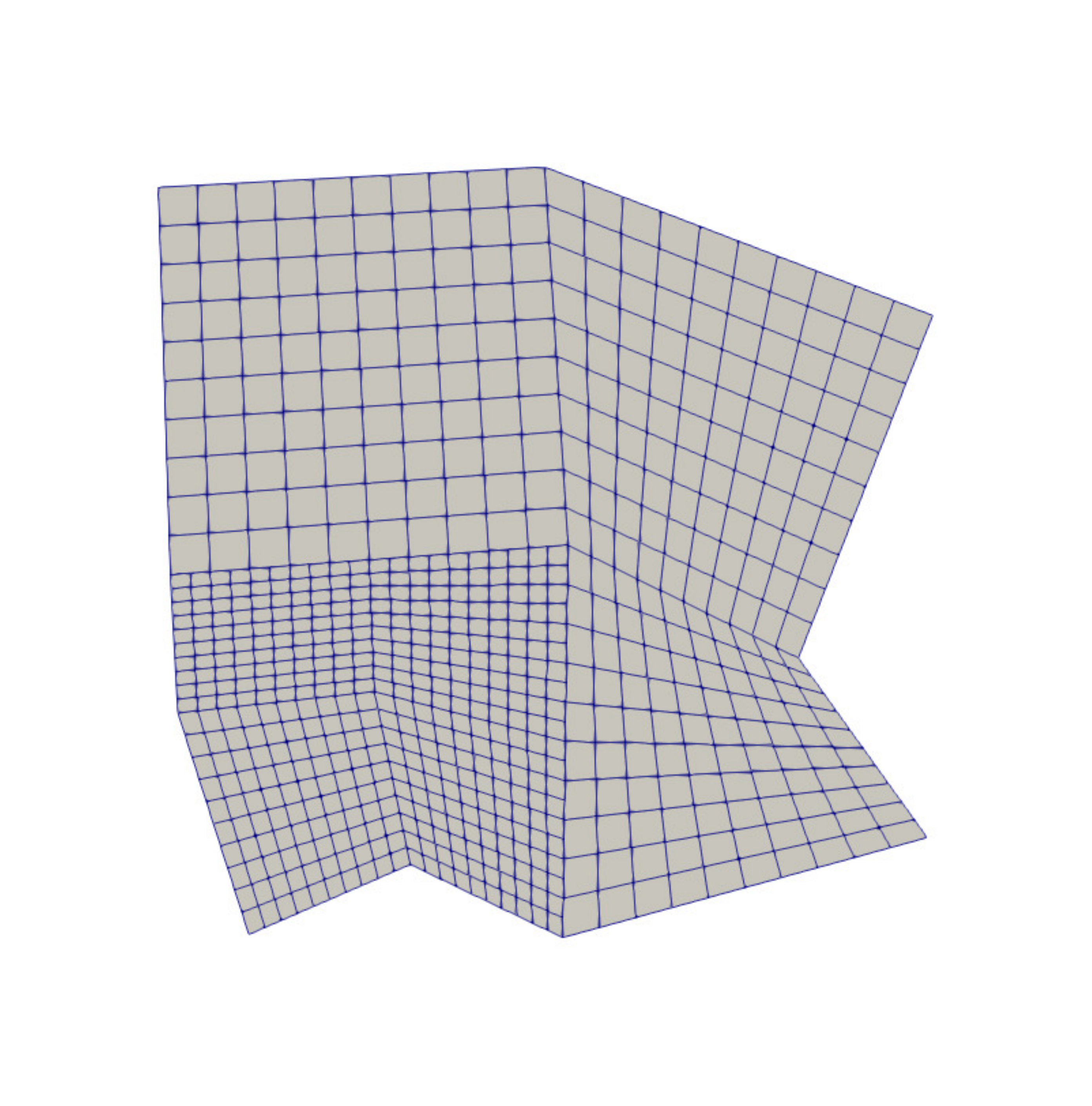}
\caption{Step 1 and 2: The bilinear mapping.}
\end{subfigure}
\begin{subfigure}{.49\linewidth}
\includegraphics[width=\linewidth]{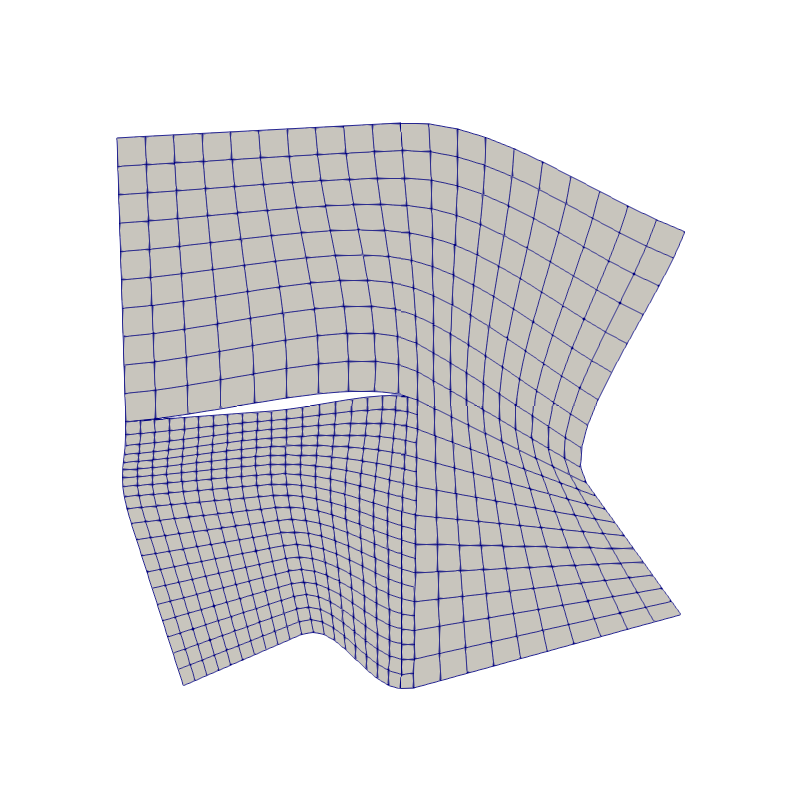}
\caption{Step 3: Averaging derivatives.}
\end{subfigure}
\begin{subfigure}{.49\linewidth}
\includegraphics[width=\linewidth]{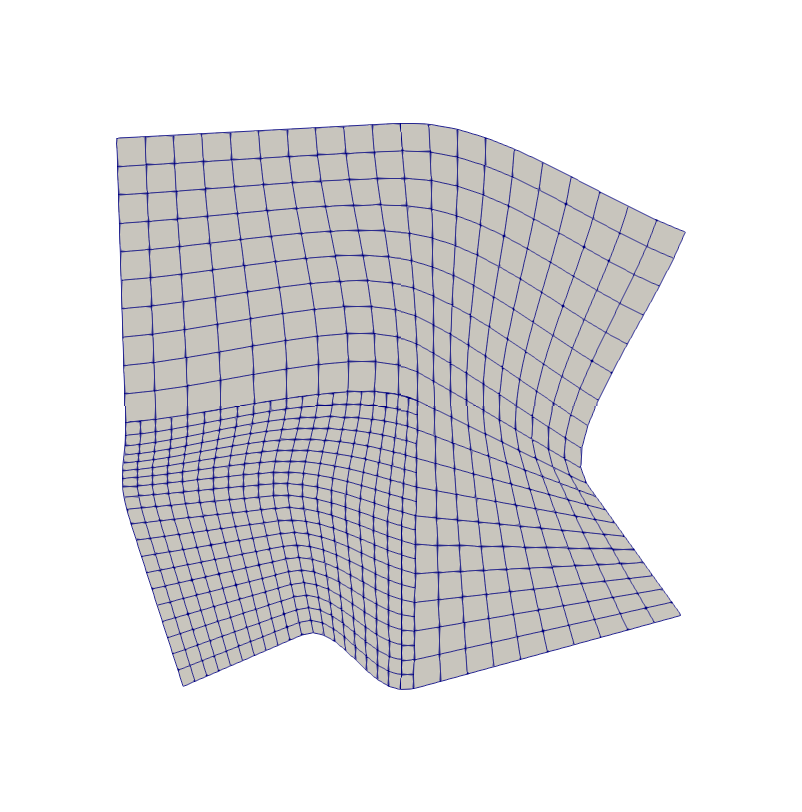}
\caption{Step 5: Computing hanging vertices.}
\end{subfigure}
\begin{subfigure}{.49\linewidth}
\includegraphics[width=\linewidth]{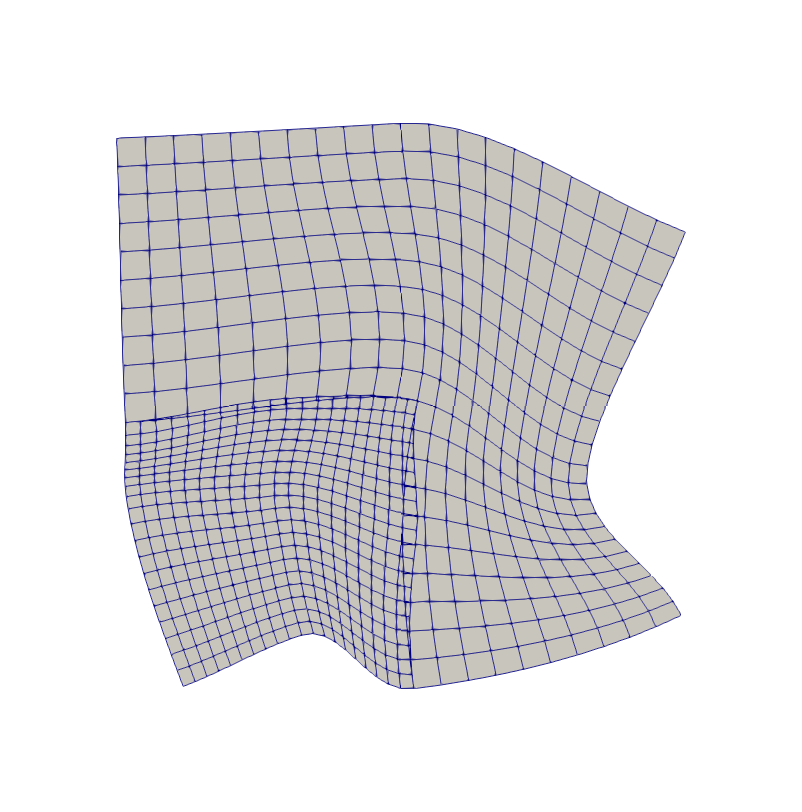}
\caption{Orthogonalization with bicubic polynomials.}
\end{subfigure}
\caption{The algorithm as performed by the \texttt{deal.II} library. Subfigures (A) to (C) illustrate the mesh obtained after steps 2, 3, and 5 without orthogonalization, respectively. Subfigure (D) shows the attempt of orthogonalization without enrichment, which is bound to fail.}
\label{fig:2D}
\end{figure}

The steps of the algorithm as performed by the \texttt{deal.II} finite
element library are illustrated in Figure~\ref{fig:2D}
(A--C). Subfigure (D) there illustrates, that it is impossible to
satisfy both orthogonality constraints and hanging node constraints
using cubic tensor product polynomials.  There are not enough degrees
of freedoms on the coarse cell left to orthogonalize the gradients at
hanging nodes. This problem can be solved by enriching the ansatz
space by ansatz functions that have zero value and unit derivative
both in $x$- and $y$-direction at all mid points of faces.

\begin{figure}
\begin{subfigure}{.49\linewidth}
\includegraphics[width=\linewidth]{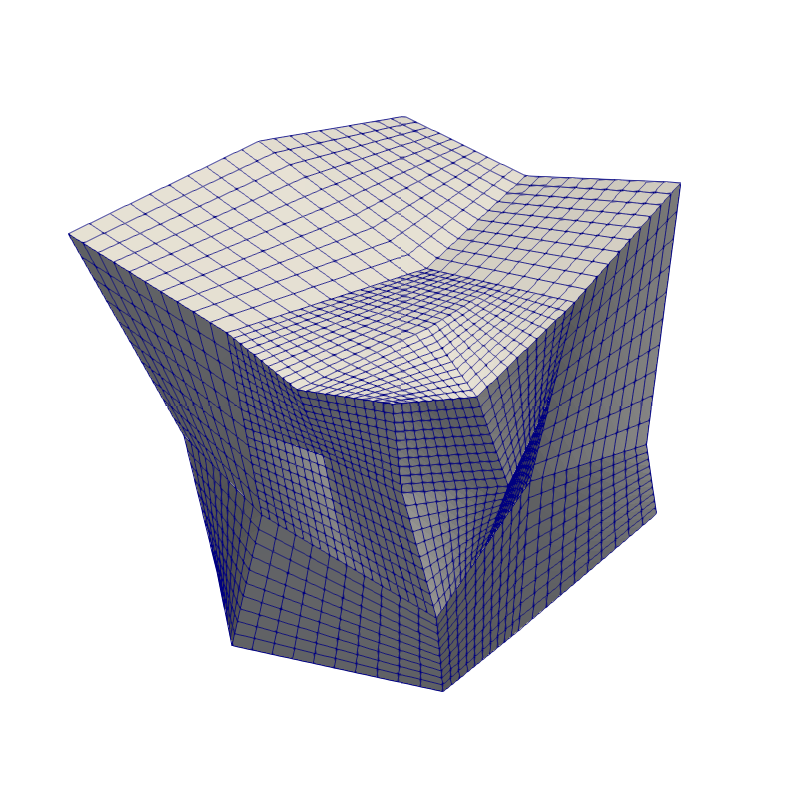}
\caption{Trilinear continuous geometry.}
\label{fig:3d-initial}
\end{subfigure}
\begin{subfigure}{.49\linewidth}
\includegraphics[width=\linewidth]{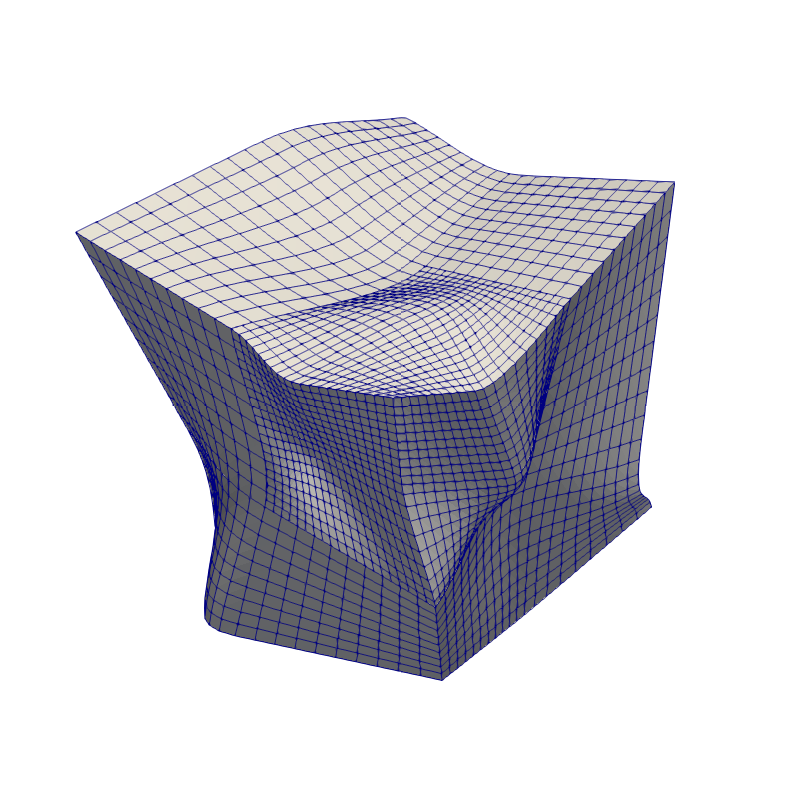}
\caption{Tricubic $C^1$-geometry.}
\label{fig:3d-mapped}
\end{subfigure}
\caption{The result of the algorithm in three dimensions computed by the \texttt{deal.II} library.}
\label{fig:3D}
\end{figure}

While the presentation above often used the two-dimensional case as an
example, the approach does not depend on the dimension. An example can
be seen in Figure~\ref{fig:3D}.

The last step of this algorithm could be deferred to the point, where
a finite element computation is executed on the mesh and coordinates
are needed on the smaller cells. In this case, the necessary values
are computed on the fly and no data is attached to the hanging
vertices. If the data is precomputed as described above, and as
implemented in \texttt{deal.II}, the higher order degrees of freedom
in these vertices may be nonzero due to interpolation and must be
stored as well.

After the degrees of freedom of the mapping in each vertex have been
computed using this algorithm, they can be used locally to determine
the geometry of each cell. On a plain, shape regular and conforming
mesh, this can be done in a straightforward way. If on the other hand
the mappings are modified for hanging vertices or anisotropic cells,
derivative degrees of freedom must be rescaled.

Rescaling on shape-regular mesh hierarchies with hanging vertices can
be achieved without storing additional data by choosing for each edge
$h_E = 2^{-\ell}$, where $\ell$ is the refinement level of the cells
adjacent to $E$. This level is unique, since hanging vertices are
ignored in the computation. Furthermore, this length is characteristic
to the cells themselves, such that the ``local'' coefficients in the
expansion~\eqref{eq:4} are obtained from the ``global'' coefficients
computed by the algorithm by multiplying by $2^{-\ell o}$, where $o$
is the number of derivatives related to the node functional.

The same principle can be applied to anisotropic cells in boundary
layers, but here the anisotropic rescaling must be stored for each
cell. Both methods can be combined.

\section*{Conclusions}

We have constructed a global $C^1$-mapping based on tensor products of cubic Hermite interpolation for quadrilateral and hexahedral meshes. Such a mapping is possible under the condition that all vertices are topologically regular. The construction easily extends to hanging vertices due to adaptive mesh refinement and to anisotropic cells in boundary layers. The algorithm can be implemented in two sweeps over the mesh using only local information in each of them.

\bibliographystyle{abbrv}
\bibliography{c1mapping}

\begin{thebibliography}{1}

\bibitem{dealII90}
G.~Alzetta, D.~Arndt, W.~Bangerth, V.~Boddu, B.~Brands, D.~Davydov,
  R.~Gassmoeller, T.~Heister, L.~Heltai, K.~Kormann, M.~Kronbichler, M.~Maier,
  J.-P. Pelteret, B.~Turcksin, and D.~Wells.
\newblock The \texttt{deal.II} library, version 9.0.
\newblock {\em Journal of Numerical Mathematics}, 2018.

\bibitem{AustinManteuffelMcCormick04}
T.~M. Austin, T.~A. Manteuffel, and S.~McCormick.
\newblock A robust multilevel approach for minimizing {$\mathbf H({\rm
  div})$}-dominated functionals in an {$\mathbf H\sp 1$}-conforming finite
  element space.
\newblock {\em Numer. Linear Algebra Appl.}, 11(2-3):115--140, 2004.

\bibitem{BangerthHartmannKanschat2007}
W.~Bangerth, R.~Hartmann, and G.~Kanschat.
\newblock {deal.II} -- a general purpose object oriented finite element
  library.
\newblock {\em ACM Trans. Math. Softw.}, 33(4):24/1--24/27, 2007.

\bibitem{BognerFoxSchmit1965}
F.~K. Bogner, R.~L. Fox, and L.~A. Schmit.
\newblock The generation of interelement compatible stiffness and mass matrices
  by the use of interpolation formulæ.
\newblock In {\em Proc. Conf. Matrix Methods in Struct. Mech.}, Wright
  Patterson AF Base, Ohio, 1965. Air Force Institute of Technology.

\bibitem{CollinSangalliTakacs2016}
A.~Collin, G.~Sangalli, and T.~Takacs.
\newblock Analysis-suitable {$G^1$} multi-patch parametrizations for {$C^1$}
  isogeometric spaces.
\newblock {\em Computer Aided Geometric Design}, 47:93--113, 2016.

\bibitem{KaplSangalliTakacs2017}
M.~Kapl, G.~Sangalli, and T.~Takacs.
\newblock The argyris isogeometric space on unstructured multi-patch planar
  domains, 2017.

\bibitem{PeteraPittman1994}
J.~Petera and J.~F.~T. Pittman.
\newblock Isoparametric hermite elements.
\newblock {\em Int. J. Numer. Meth. Engrg.}, 37:3489--3519, 1994.

\bibitem{SharmaKanschat18}
N.~Sharma and G.~Kanschat.
\newblock Convergence of an adaptive divergence-conforming discontinuous
  {G}alerkin method for the {S}tokes problem.
\newblock {\em J. Numer. Math.}, 2018.
\newblock published online.

\end{thebibliography}

\end{document}